\documentclass{article}
\usepackage{amsmath,amssymb,amsthm}
\usepackage{graphics}

\newtheorem{lemma}{Lemma}

\newcommand{\conv}{\mathop{\rm conv}\nolimits}
\newcommand{\dist}{\mathop{\rm dist}\nolimits}

\title {THE KISSING PROBLEM IN THREE DIMENSIONS}

\author {Oleg R. Musin \thanks{Institute for Math. Study of Complex Systems, Moscow State University, Moscow, Russia omusin@mail.ru}}

\begin{document}
\date{}
\maketitle

\begin{abstract}
The kissing number $k(3)$ is the maximal number of equal size nonoverlapping
spheres  in three dimensions that can touch another sphere of the same size. 
This number was the subject of a famous discussion between Isaac Newton 
and David Gregory in 1694. The first proof that $k(3)=12$ was given by Sch\"utte 
and van der Waerden only in 1953. 
In this paper we present a new solution of the Newton-Gregory problem that uses our extension of the Delsarte method. This proof relies on basic calculus and simple spherical geometry. 

\end{abstract}

\medskip

{\bf Keywords:} Kissing numbers, thirteen spheres problem, Newton-Gregory problem, Legendre polynomials, Delsarte's method

\section {Introduction}

The {\it kissing number} $k(d)$ is the highest number of equal nonoverlapping spheres in ${\bf R}^d$ that can touch another sphere of the same size. In three dimensions the kissing number problem is asking how many white billiard balls can 
{\em kiss} (touch) a black ball. 

The most symmetrical configuration, 12 billiard balls around another, is achieved if the 12 balls are placed at positions corresponding to the vertices of a regular icosahedron concentric with the central ball. However, these 12 outer balls do not kiss each other and may all be moved freely. So perhaps if you moved all of them to one side a 13th ball would possibly fit in?       

This problem was the subject of a famous discussion between Isaac Newton 
and David Gregory in 1694 (May 4, 1694; see interesting article \cite{Sz} for details of this discussion). Most reports say that Newton believed the answer was 12 balls, while Gregory 
thought that 13 might be possible. However, Casselman \cite{Cas} found some puzzling features in this story.
 
This problem is often called the {\it thirteen spheres problem}. Hoppe \cite{Hop} thought he had solved the problem in 1874. But there was a mistake - an analysis of this mistake was published by Hales in 1994 \cite{Hales}  
(see also \cite{Sz1}). 
Finally this problem was solved by Sch\"utte and van der Waerden in 1953 \cite{SvdW2}. A subsequent two-page sketch of an elegant proof was given  by Leech \cite{Lee} in 1956. Most people agree that Leech's proof is correct, but there are gaps in his exposition, many involving sophisticated spherical trigonometry. (Leech's proof was presented in the first edition  of the well-known book by Aigner and Ziegler \cite{AZ}, the authors removed this chapter from the second edition because 
a complete proof would have had to include so much spherical trigonometry.)
The thirteen spheres problem continues to be of interest, and new proofs have been published in the last few years by Hsiang \cite{Hs}, Maehara \cite{Ma},  B\"or\"oczky \cite{Bor} and Anstreicher \cite{Ans}.

The main progress in the kissing number problem in high dimensions was at the end of the 1970s. Levenshtein \cite{Lev2}, and independently, Odlyzko and Sloane \cite{OdS} (= 
\cite[Chap.13]{CS}) using Delsarte's method in 1979 proved that
$k(8)=240$ and $k(24)=196560$. This proof is surprisingly short, clean, and technically easier 
than all proofs in three dimensions. However, $d=8, 24$ are the only dimensions in which this method gives a precise result. For other cases (for instance, $d=3, 4$) the upper bounds exceed the lower. 

We found an extension of the Delsarte method in 2003 \cite{Mus}(see details in \cite{Mus2}) that allowed us to prove the bound $k(4)<25$, i.e. $k(4)=24$. This extension also yields  a proof $\; k(3)<13.$

The first version of these proofs was relatively short, but used a numerical solution of some nonconvex constrained optimization problems. Later  \cite{Mus2} these calculations were reduced to calculations of roots of polynomials in one variable. 

In this paper we present a new proof of the Newton-Gregory problem. This proof needs just basic calculus and simple spherical geometry.

\section {$k(3)=12$}

Let us recall the definition of {\em Legendre polynomials} $P_k(t)$ by the recurrence formula: 

$$P_0=1,\;\; P_1=t,\;  P_2=\frac{3}{2}t^2-\frac{1}{2}, \ldots,\; P_k=\frac {2k-1}{k}\,t\,P_{k-1}-\frac{k-1}{k}\,P_{k-2}$$
or equivalently
$$P_k(t)=\frac{1}{2^k\,k!}\,\frac{d^k}{dt^k}(t^2-1)^k \quad \mbox{ (Rodrigues' formula)}.$$


\begin {lemma}
Let $X = \{x_1, x_2,\ldots, x_n\}$ be any finite subset of the unit sphere ${\bf S}^2$ in ${\bf R}^3$.
By $\phi_{i,j}=\dist(x_i,x_j)$ we denote the spherical (angular) distance between  $x_i$ and  $x_j.$ Then
$$\sum\limits_{i=1}^n \sum\limits_{j=1}^n P_k(\cos(\phi_{i,j})) \geqslant 0.    $$
\end {lemma}

This lemma easily follows from Schoenberg's theorem \cite{Scho} for Gegenbauer (ultraspherical)  polynomials $G_k^{(d)}$. (Note that $P_k=G_k^{(3)}$.) For completeness we give a proof of Lemma 1 in the Appendix.

Let 
$$f(t) = \frac{2431}{80}t^9 - \frac{1287}{20}t^7 + \frac{18333}{400}t^5 + \frac{343}{40}t^4 - \frac{83}{10}t^3 - \frac{213}{100}t^2+\frac{t}{10} - \frac{1}{200}. $$

\medskip

\noindent{\bf Remark.}
This polynomial of degree 9 is satisfying the assumptions of the extended Delsarte's method \cite{Mus, Mus2}.
An algorithm for calculating suitable polynomials is presented in the Appendix of \cite{Mus2}.

\begin {lemma} Suppose $X = \{x_1, x_2,\ldots, x_n\} \subset {\bf S}^2$. Then
$$S(X):=\sum\limits_{i=1}^n \sum\limits_{j=1}^n f(\cos(\phi_{i,j})) \geqslant n^2.$$
\end {lemma}
\begin{proof}
The expansion of $f$ in terms of $P_k$ is
$$f = \sum\limits_{k=0}^9 {c_kP_k} = P_0 + \frac{8}{5}P_1 + \frac{87}{25}P_2 + \frac{33}{20}P_3 + \frac{49}{25}P_4 + \frac{1}{10}P_5 + \frac{8}{25}P_9.$$
We have $c_0=1,\; c_k \geqslant 0,\; k=1,2,\ldots, 9.\; $ Using Lemma 1 we get
$$ S(X)=\sum\limits_{k=0}^9 c_k \sum\limits_{i=1}^n \sum\limits_{j=1}^n P_k(\cos(\phi_{i,j}))\geqslant  
\sum\limits_{i=1}^n\sum\limits_{j=1}^n c_0P_0 =   n^2.$$
\end{proof}

If $n$ unit spheres kiss the unit sphere in ${\bf R}^3$, then the set of kissing points 
is an arrangement on the central sphere such that the (Euclidean) distance between any two points is at least 1. So the kissing number problem can be stated in another way: How many points can be placed on the surface of ${\bf S}^2$ so that the angular separation between any two points is at least $60^\circ$? 

\begin {lemma} Suppose $X = \{x_1, x_2,\ldots, x_n\}$ is a subset of $\, {\bf S}^2$ such that the angular separation 
$\phi_{i,j}$ between any two distinct points $x_i, x_j$ is at least $60^{\circ}$.
Then
$$S(X)=\sum\limits_{i=1}^n \sum\limits_{j=1}^n f(\cos(\phi_{i,j})) < 13n.$$
\end {lemma}

We give a proof of Lemma 3  in the next section.

\medskip

\medskip

\noindent{\bf Theorem.} $k(3)=12.$

\medskip

\begin{proof} Suppose $X$ is a kissing arrangement on ${\bf S}^2$ with $n=k(3)$. Then $X$  satisfies the assumptions in Lemmas 2 and 3. Therefore,
$n^2 \leqslant S(X) < 13n$. From this $n<13$ follows, i.e. $n\leqslant 12.$ From the other side we have $k(3)\geqslant 12$, showing that $n=k(3)=12.$ 
\end{proof}

\section {Proof of Lemma 3.}
We need  one fact from spherical trigonometry, namely the {\em law of cosines}: 
$$\cos{\phi} = \cos{\theta_1}\cos{\theta_2}+\sin{\theta_1}\sin{\theta_2}\cos\varphi,$$ 
for a spherical triangle $ABC$ with sides of angular lengths  $\; \theta_1,\, \theta_2,\, \phi\; $ and  $\; \angle BAC=\varphi$ (Fig. 1). For $\varphi=90^{\circ}$, this reduces 
to the spherical Pythagorean theorem: $\cos{\phi} = \cos{\theta_1}\cos{\theta_2}$.

\begin{proof} {\bf 1.} 
The polynomial $f(t)$ satisfies the following properties (see Fig.2):

\noindent $(i)\; f(t)$ is a monotone decreasing function on the interval $[-1,-t_0];$
 
\noindent $(ii)\; f(t) < 0\;\;$ for $\; t\in (-t_0,1/2]; \\$ 
where $\; f(-t_0)=0, \; t_0 \approx  0.5907$.

These properties hold because $f(t)$ has only one root $-t_0$ on $[-1, 1/2]$, and there are no zeros of the derivative $f'(t)$ (eighth degree polynomial) on $[-1,-t_0].$
\begin{center}
\begin{picture}(320,140)(-230,-70)
\thinlines
\put(-90,-54){\line(0,1){108}}
\put(90,-54){\line(0,1){108}}
\put(-90,-54){\line(1,0){180}}
\put(-90,54){\line(1,0){180}}
\put(-90,-40){\line(1,0){180}}

\thicklines
\qbezier (-90,36)(-88,38)(-82,18)
\qbezier (-82,18)(-80,12)(-78,5)
\qbezier (-78,5)(-76,-1)(-74,-7)
\qbezier (-74,-7)(-72,-13)(-70,-18)
\qbezier (-70,-18)(-66,-27) (-62,-33)

\qbezier (-62,-33)(-60,-35)(-58,-37)
\qbezier (-58,-37)(-56,-38)(-54,-39)
\qbezier (-54,-39)(-52,-40)(-50,-41)
\qbezier (-50,-41)(-46,-42)(-35,-42)
\qbezier (-35,-42)(-30,-41)(-24,-41)
\qbezier (-24,-41)(-3,-42)(10,-41) 
\qbezier (10,-41)(16,-42)(20,-43)
\qbezier (20,-43)(26,-44)(30,-45)
\qbezier (30,-45)(35,-46)(40,-48) 
\qbezier (40,-48)(45,-49)(50,-48)
\qbezier (50,-48)(55,-46)(60,-40)
\qbezier (60,-40)(66,-34)(70,-24)
\qbezier (70,-24)(75,-10)(80,6)
\qbezier (80,6)(86,24)(89,46)

\thinlines
\multiput (-80,-54)(10,0){17}%
{\line(0,1){2}}
\put(-97,-63){$-1$}
\put(-55,-63){$-0.5$}
\put(8,-63){$0$}
\put(53,-63){$0.5$}
\put(82,-63){$0.8$}
\put(-98,-42){$0$}
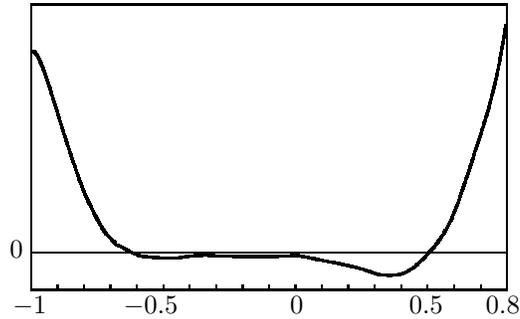
\put(-83,-80){Fig. 2. The graph of the function $f(t)$}

\thicklines
\put(-230,-54){\line(1,2){40}}
\put(-150,-54){\line(-1,2){40}}
\put(-230,-54){\line(1,0){80}}
\thinlines
\qbezier (-185,16)(-190,14)(-195,16)
\put(-194,8){$\varphi$}
\put(-193,30){$A$}
\put(-242,-56){$B$}
\put(-148,-56){$C$}
\put(-221,-14){$\theta_1$}
\put(-168,-14){$\theta_2$}
\put(-191,-50){$\phi$}

\put(-200,-80){Fig. 1}

\end{picture}
\end{center}

$$\mbox{Let} \; S_i(X):=\sum\limits_{j=1}^n f(\cos(\phi_{i,j})), \; \mbox{then} \; S(X)=\sum\limits_{i=1}^n S_i(X). \; \,  \mbox{From this it follows} $$
that if $ \; S_i(X)<13 \; $ for $ \; i=1, 2, \ldots, n, \; $ then $\; S(X)<13n$. 

We obviously have $\phi_{i,i}=0$, so $f(\cos{\phi_{i,i}})=f(1)$. Note that our assumption on  $X$ ($\phi_{i,j} \geqslant 60^{\circ}, \,  i \ne j$) yields
$\cos{\phi_{i,j}} \leqslant 1/2.$ Therefore, $\cos{\phi_{i,j}}$ lies in the interval [-1,1/2].
By $(ii)$ we have $f(\cos{\phi_{i,j}}) \leqslant 0$ whenever $\cos{\phi_{i,j}} \in [-t_0,1/2]$. Let $J(i):=\{j:\cos{\phi_{i,j}}\in [-1,-t_0)\}$. We obtain  
\begin{equation}
\label{one}
S_i(X) \leqslant T_i(X):=f(1)+\sum\limits_{j \in J(i)} f(\cos{\phi_{i,j}}). 
\end{equation}

Let $\theta_0=\arccos{t_0} \approx 53.794^\circ.$ Then $j\in J(i)$ iff 
$\phi_{i,j}>180^\circ-\theta_0$, i.e.
$\theta_j<\theta_0$, where $\theta_j=180^\circ-\phi_{i,j}.$ In other words all $x_{i,j},\; j\in J(i),$ lie inside the spherical cap of center $e_0$ and radius $\theta_0$, where $e_0=-x_i$ is the antipodal point to $x_i$.

\medskip

{\bf 2.} Let us consider on ${\bf S}^2$ points $e_0,y_1,\ldots,y_m$  such that
$$
\phi_{i,j}:=\dist(y_i,y_j)\geqslant 60^\circ, \;  \forall \; i\neq j, \quad
\theta_i:=\dist(e_0,y_i)< \theta_0 \; \mbox{ for } \; 1\leqslant i\leqslant m.
\eqno (2)$$

Denote by $\mu $ the highest value of $m$ such that the constraints in $(2)$ allow a nonempty set of points $y_1,\ldots,y_m.$ 

Suppose that $\; 0\leqslant m\leqslant\mu\;$ and $Y=\{y_1,\ldots, y_m\}$ satisfies $(2)$. Let 
$$H(Y)=H(y_1,\ldots,y_m):=f(1)+f(-\cos{\theta_1})+\ldots+f(-\cos{\theta_m}), $$ 
$$h_m:=\sup\limits_Y{\{H(Y)\}},\quad
h_{max}:=\max{\{h_0,h_1,\ldots,h_\mu\}}.$$

It is clear that $\; T_i(X)\leqslant h_m$, where $m=|J(i)|$.  From $(1)$ it follows that $ S_i(X) \leqslant h_m.$ Thus, if we prove that $h_{max} < 13$, then we prove Lemma 3. 

\medskip

{\bf 3.} Now we prove that $\mu \leqslant 4.$

\noindent Suppose $Y=\{y_1,\ldots,y_m\} \subset {\bf S}^2$ satisfies $(2)$. By symmetry we may assume that $e_0$ is the North pole and $y_i$ has polar coordinates $(\theta_i,\varphi_i)$. Then from the law of cosines we have:
$$\cos{\phi_{i,j}} = \cos{\theta_i}\cos{\theta_j}+\sin{\theta_i}\sin{\theta_j}\cos(\varphi_i-\varphi_j).$$
Note that $\theta_i>0$ for $m\geqslant 2$. Conversely, $y_i=e_0, \; \theta_j=\phi_{i,j}\geqslant 60^\circ>\theta_0$, a 
contradiction. From $(2)$ we have $\cos{\phi_{i,j}}\leqslant 1/2$, then
$$\cos(\varphi_i-\varphi_j)\leqslant \frac{1/2-\cos{\theta_i}\cos{\theta_j}}{\sin{\theta_i}\sin{\theta_j}}.\eqno (3)$$
Let  $$Q(\alpha,\beta):=\frac{1/2-\cos{\alpha}\cos{\beta}}{\sin{\alpha}\sin{\beta}},$$   then     
$$Q'_\alpha(\alpha,\beta)=
\frac{\partial{Q(\alpha,\beta)}}{\partial{\alpha}}=\frac{2\cos{\beta}-\cos{\alpha}}{2\sin^2{\alpha}\sin{\beta}}.$$ 
From this it follows that if $\; 0<\alpha, \beta\leqslant \theta_0$, then $ \cos{\beta}> 1/2$ (because $\theta_0<60^\circ$); so then $Q'_\alpha(\alpha,\beta)> 0,$ and $Q(\alpha,\beta)\leqslant 
Q(\theta_0,\beta)=Q(\beta,\theta_0)\leqslant Q(\theta_0,\theta_0).$ Therefore,
$$\frac{1/2-\cos{\theta_i}\cos{\theta_j}}{\sin{\theta_i}\sin{\theta_j}} \leqslant
\frac{1/2-\cos^2{\theta_0}}{\sin^2{\theta_0}}=\frac{1/2-t_0^2}{1-t_0^2}.$$

Combining this inequality and (3), we get $$\cos(\varphi_i-\varphi_j)\leqslant 
\frac{1/2-t_0^2}{1-t_0^2}.$$ 
Note that $\arccos((1/2-t_0^2)/(1-t_0^2)) \approx 76.582^\circ>72^\circ$. This implies that $m\leqslant 4$ because no more than four points can lie in a circle with  minimum angular separation between any two points greater than $72^\circ$.

\medskip

{\bf 4.} Now we have to prove that $h_{max}=\max{\{h_0,h_1,h_2,h_3,h_4\}} < 13.\\$
We obviously have $h_0=f(1)=10.11<13$.

From $(i)$ follows that $f(-\cos{\theta})$ is a monotone decreasing function in $\theta$ on 
$[0,\theta_0].$  
Then for $m=1:\; H(y_1)=f(1)+f(-\cos{\theta_1})$ attains its maximum at $\theta_1=0.$ So then
$$h_1=f(1)+f(-1)=12.88<13.$$ 

{\bf 5.}
Let us consider for $m=2,3,4$ an  arrangement
$\{e_0,y_1,\ldots,y_m\}$ in ${\bf S}^2$ that gives $H(Y)=h_m$. Here $y_i\neq e_0$ (see {\bf 3}).
Note that in this arrangement, points $y_k$  cannot be shifted towards $e_0$ because in this case $H(Y)$ increases. 

For $m=2$ this yields $e_0\in y_1y_2,$ and  $\dist(y_1,y_2)= 60^\circ.$ If  $e_0\notin y_1y_2,$
then the whole arc $y_1y_2$ can be shifted towards $e_0$. If $\dist(y_1,y_2)> 60^\circ,$ then
$y_1$ (and $y_2$) can be shifted towards $e_0.$

For $m=3$ we prove that $\Delta_3=y_1y_2y_3$ is a spherical regular triangle with edge length $60^\circ$. As above, $e_0\in \Delta_3$, otherwise the whole triangle can be shifted towards $e_0$. Suppose $\dist(y_1,y_i)>60^\circ,\; i=2,3,$ then $\dist(y_1,e_0)$ can be decreased. From this follows that for any $y_i$ at least one of the distances $\{\dist(y_i,y_j)\}$ is equal to $60^\circ$. Therefore, at least two sides of 
$\Delta_3$ (say $y_1y_2$ and $y_1y_3$) have length $60^\circ$. Also $\dist(y_2,y_3)=60^\circ$, conversely  $y_3$ (or $y_2$, if $e_0 \in y_1y_3$) can be rotated about $y_1$ by a small angle towards $e_0$ (Fig.3).

For $m=4$ we first prove that $\Delta_4:=\conv{Y}$ (the spherical convex hull of $Y$) is a convex quadrilateral. 
Conversely, we may assume that $y_4\in y_1y_2y_3$. 

The great circle through $y_4$ that is orthogonal to the arc $e_0y_4$ divides ${\bf S}^2$ into two hemispheres: $H_1$ and $H_2$.  Suppose 
$e_0\in H_1$,   then at least one $y_i$ (say $y_3$) belongs to $H_2$ (Fig.4).   
So the angle $\; \angle {e_0y_4y_3}\; $ greater  than $90^\circ$, then (again from the law of cosines)  $\quad \dist(y_3,e_0)>\dist(y_3,y_4). \quad $ Thus,
$\\ \theta_3=\dist(y_3,e_0)>\dist(y_3,y_4)\geqslant 60^\circ>\theta_0 \quad - \; \, $ a contradiction.

Arguing as for $m=3$ it is easy to prove that for any vertex $y_i$ there are at least two vertices $y_j$ at the distance $60^\circ$ from $y_i$. Note that the diagonals of  $\Delta_4$ cannot be both of lengths 
$60^\circ.$ Conversely, at least one side of $\Delta_4$ is of length less than $60^\circ.$ Thus, $\Delta_4$ is a spherical equilateral quadrangle (rhomb) with edge length $60^\circ$. 

\medskip

{\bf 6.} Now we introduce the function $F_1(\psi),$\footnote{For given $\psi$, the value $F_1(\psi)$
can be found as the maximum of the  9th degree polynomial  $\Omega(s)=\widetilde F_1(\theta,\psi), \, s=\cos{(\theta-\psi/2)},$ on the interval $[\cos(\theta_0-\psi/2),1].$} 
where $\psi\in [60^\circ,2\theta_0]$:
$$F_1(\psi):=\max\limits_{\psi/2\leqslant \theta \leqslant \theta_0}
\{\widetilde F_1(\theta,\psi)\}, \quad 
\widetilde F_1(\theta,\psi)=f(-\cos{\theta})+f(-\cos(\psi-\theta)).$$
So if $\, \dist(y_i,y_j)=\psi,\;$ then
$$f(-\cos{\theta_i})+f(-\cos{\theta_j})\leqslant F_1(\psi).\eqno (4)$$ 

Therefore,

$$H(y_1,y_2)\leqslant h_2=f(1)+F_1(60^\circ)\approx 12.8749 < 13.$$


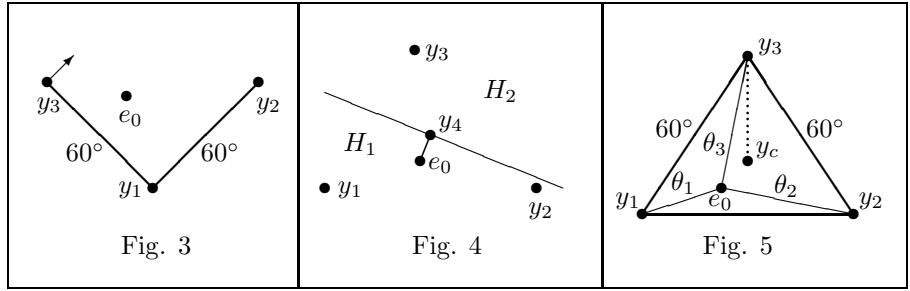
\begin{center}
\begin{picture}(320,110)(0,180)

\thicklines
\put(45,220){\line(-1,1){40}}
\put(45,220){\line(1,1){40}}

\thinlines

\put(5,260){\vector(1,1){10}}

%

\put(45,220){\circle*{4}}
\put(5,260){\circle*{4}}
\put(85,260){\circle*{4}}
\put(35,255){\circle*{4}}
\put(32,245){$e_0$}
\put(32,218){$y_1$}
\put(1,250){$y_3$}
\put(84,250){$y_2$}
\put(12,230){$60^\circ$}
\put(63,230){$60^\circ$}

\put(33,195){Fig. 3}

\thicklines
\put(230,210){\line(2,3){40}}
\put(310,210){\line(-2,3){40}}
\put(230,210){\line(1,0){80}}

\thinlines
\multiput(270,233.5)(0,2.5){15}%
{\circle*{1}}
\put(260,220){\line(1,5){10}}
\put(260,220){\line(-3,-1){30}}
\put(260,220){\line(5,-1){50}}

\put(270,230){\circle*{4}}
\put(230,210){\circle*{4}}
\put(310,210){\circle*{4}}
\put(270,270){\circle*{4}}
\put(260,220){\circle*{4}}

\put(220,214){$y_1$}
\put(312,214){$y_2$}
\put(274,272){$y_3$}
\put(273,233){$y_c$}
\put(255,213){$e_0$}
\put(241,219){$\theta_1$}
\put(280,218){$\theta_2$}
\put(253,234){$\theta_3$}
\put(235,239){$60^\circ$}
\put(292,239){$60^\circ$}
\put(253,195){Fig. 5}

\put(-10,290){\line(1,0){340}}
\put(-10,182){\line(1,0){340}}

\put(-10,182){\line(0,1){108}}
\put(330,182){\line(0,1){108}}

\put(100,182){\line(0,1){108}}
\put(215,182){\line(0,1){108}}


\put(150,240){\circle*{4}}
\put(110,220){\circle*{4}}
\put(190,220){\circle*{4}}
\put(146,230){\circle*{4}}

\put(115,218){$y_1$}
\put(187,210){$y_2$}
\put(148,270){$y_3$}
\put(144,272){\circle*{4}}

\put(153,243){$y_4$}
\put(149,227){$e_0$}

\put(110,256){\line(5,-2){90}}
\put(150,240){\line(-2,-5){4}}
\put(170,254){$H_2$}
\put(117,234){$H_1$}
\put(143,195){Fig. 4}




\end{picture}
\end{center}

{\bf 7.} Let $m=4,$ 
$\; d_1=\dist(y_1,y_3), \; d_2=\dist(y_2,y_4),\; $ Since  $\; \Delta_4=y_1y_2y_3y_4$ is a spherical rhomb, we have
$\; \cos(d_1/2)\cos(d_2/2)=1/2\,$ (Pythagorean theorem, the diagonals $y_1y_3,\, y_2y_4$ of $\Delta_4$ are orthogonal). So if
$$\rho(s):=2\arccos{\frac{1}{2\cos(s/2)}},$$ 
then $$\rho(d_1)=d_2, \quad  \rho(d_2)=d_1,\quad \rho(90^\circ)=90^\circ, \quad \rho(\rho(s))=s.$$

 Suppose  $d_1\leqslant d_2$. The inequalities  $\; \theta_i\leqslant \theta_0\; $ yield  $\; d_2\leqslant 2\theta_0.\; $  Then   $$\rho(2\theta_0)\leqslant d_1 \leqslant 90^\circ\leqslant d_2\leqslant 2\theta_0.$$

Now we consider two cases: $\\1) \; \rho(2\theta_0) \leqslant d_1<77^\circ, \; $ and $\; 2) \; 77^\circ\leqslant d_1\leqslant 90^\circ.\\$ 
$1)\;$ Clearly, $F_1(\psi)$ is a monotone decreasing function in $\psi.\;$ Then $(4)$ implies
$$f(-\cos{\theta_1})+f(-\cos{\theta_3})\leqslant F_1(d_1)\leqslant F_1(\rho(2\theta_0)),$$ 
$$ f(-\cos{\theta_2})+f(-\cos{\theta_4})\leqslant F_1(d_2)=F_1(\rho(d_1)) < F_1(\rho(77^\circ)),$$
so then
$$H(Y) < f(1)+F_1(\rho(2\theta_0))+F_1(\rho(77^\circ))\approx 12.9171<13.$$
$2)$ In this case we have
$$H(Y) \leqslant f(1)+F_1(77^\circ)+F_1(90^\circ)\approx 12.9182<13.$$

Thus, $\; h_4 <13.$

\medskip




{\bf 8.} Our last step is to show that $h_3<13.$\footnote{A more detailed analysis shows $h_3\approx 12.8721, \; h_4\approx 12.4849.$}

Since $\Delta_3$ is a regular triangle, $H(Y)=f(1)+f(-\cos{\theta_1})+f(-\cos{\theta_2})+f(-\cos{\theta_3})$ is a symmetric function in the $\theta_i$, so it is sufficient to consider  the case
$\; \theta_1\leqslant \theta_2\leqslant \theta_3\leqslant \theta_0. $ 

In this case
$\; R_0\leqslant \theta_3\leqslant\theta_0,\; $ where $\; R_0=\arccos{\sqrt{2/3}}\approx  35.2644^\circ \; $  is the (spherical) circumradius of $\Delta_3$.

Let  $y_c$ be the center of $\Delta_3$. We have $\gamma:=\angle {y_1y_3y_c}=\angle {y_2y_3y_c}.$ Using the low of cosines for the triangle $y_1y_3y_c$, we get $\; \gamma=\arccos{\sqrt{2/3}}, \; $
i.e. $\; \gamma=R_0.$

Denote  the angle $\angle {e_0y_3y_c}$ by $u$. Then (see Fig.5)
$$\cos{\theta_1}=\cos{60^\circ}\cos{\theta_3}+\sin{60^\circ}\sin{\theta_3}\cos{(R_0-u)},$$
$$\cos{\theta_2}=\cos{60^\circ}\cos{\theta_3}+\sin{60^\circ}\sin{\theta_3}\cos{(R_0+u)},$$ 
where $\; 0\leqslant u\leqslant u_0:=\arccos(\cot{\theta_3}/\sqrt{3})-R_0.$ Note that if $\, u=u_0, \, $ then $\, \theta_2=\theta_3; \; \, u=0 \; $ yields $\; \theta_1=\theta_2; \; $ and  if $\; \, 0<u<u_0,\; $ then $\; \theta_1<\theta_2<\theta_3$. 

For fixed $\theta_3=\psi,\; H(y_1,y_2)$ is a polynomial of degree 9 in $s=\cos{u}$. Denote by 
$F_2(\psi)$ the maximum of this polynomial on the interval $[\cos{u_0},1]$. 

Let
$$\{\psi_1,\ldots,\psi_6\}=\{R_0,\, 38^\circ,\, 41^\circ,\, 44^\circ,\, 48^\circ,\, \theta_0\}.\;$$
It is clear that $F_2(\psi)$  is a monotone increasing function in $\psi$ on  $[R_0,\theta_0].\; $ 
From the other side, $f(-\cos{\psi})$ is a monotone decreasing function in $\psi$.
Therefore for  $\theta_3\in [\psi_i,\psi_{i+1}]$ we have
$$H(Y)=H(y_1,y_2)+f(-\cos{\theta_3}) < w_i:=F_2(\psi_{i+1})+f(-\cos{\psi_i}).$$ 
Since,
$$\; \{w_1,\ldots,w_5\}\approx\{12.9425, 12.9648, 12.9508, 12.9606, 12.9519\},$$ 
we get $\; h_3<\max\{w_i\}<13.$

Thus, $h_m <13$ for all $m$ as required.
\end{proof}

\medskip

\medskip

\medskip

\noindent {\Large\bf Appendix. Proof of Lemma 1.}

\medskip

In this proof we are using Schoenberg's original proof \cite{Scho} which is based on the addition theorem for Gegenbauer polynomials.\footnote{Pfender and Ziegler\cite{PZ} give a proof as a simple consequence of the addition theorem for spherical harmonics. This theorem is not so elementary. The addition theorem for Legendre polynomials can be proven by elementary algebraic calculations.} 
The addition theorem for Legendre polynomials was discovered by Laplace and Legendre in 1782-1785:
$$P_k(\cos{\theta_1}\cos{\theta_2}+\sin{\theta_1}\sin{\theta_2}\cos{\varphi})$$
 $$=P_k(\cos{\theta_1})P_k(\cos{\theta_2})
+ \, 2\sum\limits_{m=1}^k\,\frac{(k-m)!}{(k+m)!}\,P_k^m(\cos{\theta_1})P_k^m(\cos{\theta_2})\,\cos{m\varphi}$$
$$= \, \sum\limits_{m=0}^k c_{m,k}\,P_k^m(\cos{\theta_1})P_k^m(\cos{\theta_2})\,\cos{m\varphi},$$
where
$$P_k^m(t) = (1-t^2)^{\frac{m}{2}}\,\frac{d^m}{dt^m}P_k(t).$$
(See details in \cite{Car} and  \cite {Erd}.)

\medskip

\begin{proof}
Let $X=\{x_1, \ldots, x_n\} \subset {\bf S}^2$ and $x_i$ has spherical (polar) coordinates $(\theta_i,\varphi_i)$. Then from the law of cosines we have:
$$\cos{\phi_{i,j}}=\cos{\theta_i}\,\cos{\theta_j}+\sin{\theta_i}\sin{\theta_j}\cos{\varphi_{i,j}},\quad
\varphi_{i,j}:=\varphi_i-\varphi_j,$$
which yields
$$\sum\limits_{i,j}P_k(\cos{\phi_{i,j}})=\sum\limits_{i,j}\sum\limits_{m=0}^kc_{m,k}P_k^m(\cos{\theta_i})P_k^m(\cos{\theta_j})\cos{m\varphi_{i,j}}$$
$$ = \sum\limits_mc_{m,k}\sum\limits_{i,j}u_{m,i}u_{m,j}\cos{m\varphi_{i,j}}, \quad 
u_{m,i}=P_k^m(\cos{\theta_i}).$$ 

Let us prove that for any real $u_1, \ldots, u_n$
$$\sum_{i,j}u_iu_j\cos{m\varphi_{i,j}}\geqslant 0.$$

Pick $n$ vectors $v_1, \ldots, v_n$  in ${\bf R}^2$ with coordinates 
$v_i=(\cos{m\varphi_i}, \sin{m\varphi_i})$. If $v=u_1v_1+\ldots+u_nv_n,$ then
$$0\, \leqslant  \, ||v||^2 \, = \, \langle v,v\rangle \, = \, \sum_{i,j}u_iu_j\cos{m\varphi_{i,j}}.$$
This inequality and the inequalities $c_{m,k}>0$ complete our proof.
\end{proof}


\begin{thebibliography}{99}

\bibitem{AZ}
M. Aigner and G.M. Ziegler, Proofs from THE BOOK, Springer, 1998 (first ed.) and 2002 (second ed.)
\bibitem{Ans}
K. Anstreicher, The thirteen spheres: A new proof, Discrete and Computational Geometry, {\bf 31}(2004), 613-625.
\bibitem{Bor}
K. B\"or\"oczky, The Newton-Gregory problem revisited, Proc. Discrete Geometry, Marcel Dekker, 2003,  103-110.
\bibitem{Car}
B.C. Carlson, Special functions of applied mathematics, Academic Press, 1977.
\bibitem{Cas}
B. Casselman, The difficulties of kissing in three dimensions, Notices Amer. Math. Soc., {\bf 51}(2004), 884-885. 
\bibitem{CS}
J.H. Conway and N.J.A. Sloane, Sphere Packings, Lattices, and Groups, New York, Springer-Verlag, 1999 (Third Edition).

\bibitem{Erd}
A. Erd\'elyi, editor, Higher Transcendental Function, McGraw-Hill, NY, 3 vols, 1953, Vol. II,
Chap. XI. 
\bibitem{Hales}
T. Hales, The status of the Kepler conjecture, Mathematical Intelligencer {\bf 16}(1994), 47-58.

\bibitem{Hop}
R. Hoppe, Bemerkung der Redaction, Archiv Math. Physik (Grunet) {\bf 56} (1874), 307-312.    

\bibitem{Hs}
W.-Y. Hsiang, Least action principle of crystal formation of dense packing type and Kepler's conjecture, World Scientific, 2001. 

\bibitem{Lee}
J. Leech, The problem of the thirteen spheres, Math. Gazette {\bf 41} (1956), 22-23. 
\bibitem{Lev2}
V.I. Levenshtein, On bounds for packing in $n$-dimensional Euclidean space, Sov. Math. Dokl. 
{\bf 20}(2), 1979, 417-421.

\bibitem{Ma}
H. Maehara, Isoperimetric theorem for spherical polygons and the problem of 13 spheres, Ryukyu Math. J.,
{\bf 14} (2001), 41-57.

\bibitem{Mus}
O.R. Musin, The problem of the twenty-five spheres, Russian Math. Surveys, {\bf 58}(2003), 794-795. 

\bibitem{Mus2}
O.R. Musin, The kissing number in four dimensions, preprint, September 2003, math. MG/0309430.

\bibitem{OdS}
A.M. Odlyzko and N.J.A. Sloane, New bounds on the number of unit spheres that
can touch a unit sphere in $n$ dimensions, J. of Combinatorial Theory
A26(1979), 210-214.

\bibitem{PZ}
F. Pfender and G.M. Ziegler, Kissing numbers, sphere packings, and some unexpected proofs, Notices Amer. Math. Soc., {\bf 51}(2004), 873-883. 

\bibitem{Scho}
I.J. Schoenberg,  Positive definite functions on spheres, Duke Math. J., 
{\bf 9} (1942), 96-107.

\bibitem{SvdW2}
K. Sch\"utte and B.L. van der Waerden, Das Problem der dreizehn Kugeln, Math. Ann. {\bf 125} (1953), 325-334.
\bibitem{Sz1}
G.G. Szpiro, Kepler's conjecture, Wiley, 2002.
\bibitem{Sz}
G.G. Szpiro, Newton and the kissing problem,\\ http://plus.maths.org/issue23/features/kissing/


 \end{thebibliography}
\end{document}